# VOLTERRA PERIODIC CHAIN

## Belskiy D. V.

*Institute of Mathematics NAS of Ukraine*
*Ukraine, 01024, Kyiv, st. Tereshchenkovskaya, 3*
*e-mail: oiop120@gmail.com*

In this paper, the inverse spectral problem is applied to the integration of a periodic Volterra chain. A generalization of the Lagrange interpolation formula has been made.

**1. Introduction.** Recently, nonlinear evolutionary equations of the KdV type have attracted great interest in connection with various physical applications. Among such equations is the Volterra chain [1, 2]

$$u'_n = u_n \left( u_{n+1} - u_{n-1} \right), \ n \in \mathbb{Z}. \tag{1}$$

In the case of $u_n > 0$, this system, in particular, arises in the study of the fine structure of the spectra of Langmuir oscillations in plasma, which is why it is sometimes called the Langmuir chain [1]. System (1) can also be considered as a special case of Volterra mathematical models in biology [3, 4]. The Volterra chain has been studied in many works (see [5] and the literature cited there); in particular, the periodic Volterra chain was investigated in [6-8]. In this paper, a method for integrating the periodic Volterra chain by differentiating the spectral parameters of the discrete Hill equation [9-11] is presented, which is a continuation of the studies begun in [6] and is substantially based on the results obtained in [9-12]. As an auxiliary result, a generalization of the Lagrange interpolation formula is made. The goal of this paper is to obtain representations for solutions of the periodic system (1) within the framework of the inverse spectral problem for the discrete Hill equation.

**2. Lagrange's interpolation formula.** In the next lemma we will make a generalization of the Lagrange's interpolation formula (see [11], Appendix E), which we will need later.

**Lemma.** *For quantities* $s \in \mathbb{Z}$, $x_j \in \mathbb{C}$, $j = \overline{0,n}$, $n \in \mathbb{N}$, *such that* $x_j \neq x_m$ *when* $j \neq m$, *the equality*

$$\sum_{j=0}^{n} \frac{x_j^s}{\prod_{\substack{l=0 \\ l \neq j}}^{n} (x_j - x_l)} = \begin{cases} \dfrac{(-1)^n}{\prod_{j=0}^{n} x_j} g_{|s|-1}\left(\dfrac{1}{x_0}, \dfrac{1}{x_1}, \ldots, \dfrac{1}{x_n}\right), \ s \leq -1, \ n \in \mathbb{N}; \\ 0, \ s \geq 0, \ n \geq s+1; \\ g_k(x_0, x_1, \ldots, x_n), \ s = n+k, \ k \in \mathbb{N} \cup \{0\}, \ n \in \mathbb{N}; \end{cases}$$

*is satisfied, where* $g_l(u_0, u_1, \ldots, u_n) = \sum_{\substack{m_0 \geq 0, m_1 \geq 0, \ldots, m_n \geq 0 \\ m_0 + m_1 + \ldots + m_n = l}} u_0^{m_0} u_1^{m_1} \ldots u_n^{m_n}$, $\{m_0, m_1, \ldots, m_n, l\} \subset \mathbb{N} \cup \{0\}$; $x_j \neq 0$, $j = \overline{0,n}$, *when* $s \leq -1$.

*Proof.* The first method. We write the sum in the formulation of the lemma as follows

$$\sum_{j=0}^{n} \frac{x_j^s}{\prod_{\substack{l=0 \\ l \neq j}}^{n} (x_j - x_l)} = \sum_{j=0}^{n-1} \frac{x_j^s}{\prod_{\substack{l=0 \\ l \neq j}}^{n} (x_j - x_l)} + \frac{x_n^s}{\prod_{\substack{l=0 \\ l \neq n}}^{n} (x_n - x_l)} =$$

$$= \sum_{j=0}^{n-1} \frac{x_j^s}{\left[\prod_{\substack{l=0 \\ l \neq j}}^{n-1} (x_j - x_l)\right](x_j - x_n)} + \frac{x_n^s}{\prod_{\substack{l=0 \\ l \neq n}}^{n} (x_n - x_l)} =$$

$$= \sum_{j=0}^{n-1} \left[ \frac{x_j^{s-1}}{\prod_{\substack{l=0 \\ l \neq j}}^{n-1}(x_j - x_l)} \left(1 + \frac{x_n}{x_j - x_n}\right) \right] + \frac{x_n^s}{\prod_{\substack{l=0 \\ l \neq n}}^{n}(x_n - x_l)} =$$

$$= \sum_{j=0}^{n-1} \frac{x_j^{s-1}}{\prod_{\substack{l=0 \\ l \neq j}}^{n-1}(x_j - x_l)} + \sum_{j=0}^{n-1} \left[ \frac{x_j^{s-1}}{\prod_{\substack{l=0 \\ l \neq j}}^{n-1}(x_j - x_l)} \frac{x_n}{x_j - x_n} \right] + \frac{x_n^s}{\prod_{\substack{l=0 \\ l \neq n}}^{n}(x_n - x_l)} =$$

$$= \sum_{j=0}^{n-1} \frac{x_j^{s-1}}{\prod_{\substack{l=0 \\ l \neq j}}^{n-1}(x_j - x_l)} + x_n \sum_{j=0}^{n-1} \frac{x_j^{s-1}}{\prod_{\substack{l=0 \\ l \neq j}}^{n}(x_j - x_l)} + \frac{x_n^{s-1}}{\prod_{\substack{l=0 \\ l \neq n}}^{n}(x_n - x_l)} \cdot x_n =$$

$$= \sum_{j=0}^{n-1} \frac{x_j^{s-1}}{\prod_{\substack{l=0 \\ l \neq j}}^{n-1}(x_j - x_l)} + x_n \sum_{j=0}^{n} \frac{x_j^{s-1}}{\prod_{\substack{l=0 \\ l \neq j}}^{n}(x_j - x_l)}.$$

Denoting

$$f(s,n) \stackrel{df}{=} \sum_{j=0}^{n} \frac{x_j^s}{\prod_{\substack{l=0 \\ l \neq j}}^{n}(x_j - x_l)},$$

we write the last relationship as

$$f(s,n) = f(s-1, n-1) + f(s-1, n) x_n. \qquad (2)$$

For $s \leq -1$ the equality

$$f(s,1) = \frac{x_0^s}{x_0 - x_1} + \frac{x_1^s}{x_1 - x_0} = -\frac{1}{x_0^{|s|} x_1^{|s|}} \frac{x_0^{|s|} - x_1^{|s|}}{x_0 - x_1} =$$

$$= -\frac{1}{x_0 x_1^{|s|}} \left[ 1 + \frac{x_1}{x_0} + \left(\frac{x_1}{x_0}\right)^2 + \ldots + \left(\frac{x_1}{x_0}\right)^{|s|-1} \right] = -\frac{1}{x_0 x_1} g_{|s|-1}\left(\frac{1}{x_0}, \frac{1}{x_1}\right) \qquad (3)$$

holds. In [11], on page 188, it is shown that $f(0,n) = 0$ for all $n \in \mathbb{N}$. Therefore, from (2) for $n \geq 2$ we obtain

$$0 = f(0,n) = f(-1, n-1) + f(-1, n) x_n,$$

$$f(-1, n) = -\frac{1}{x_n} f(-1, n-1). \qquad (4)$$

In this case, from (3) follows the equality

$$f(-1, 1) = -\frac{1}{x_0 x_1}. \qquad (5)$$

Using mathematical induction from relations (4) and (5) we obtain the identity

$$f(-1, n) = \frac{(-1)^n}{\prod_{j=0}^{n} x_j}, \; n \in \mathbb{N}. \qquad (6)$$

Let us assume that the equality

$$f(s,n) = \frac{(-1)^n}{\prod_{j=0}^{n} x_j} g_{|s|-1}\left(\frac{1}{x_0}, \frac{1}{x_1}, ..., \frac{1}{x_n}\right) \tag{7}$$

is satisfied for some $s \leq -1$ for all $n \geq 1$. In this case, from (3) we obtain the relation

$$f(s-1,1) = -\frac{1}{x_0 x_1} g_{|s|}\left(\frac{1}{x_0}, \frac{1}{x_1}\right). \tag{8}$$

Now let us assume that for some $n \geq 2$ the equality

$$f(s-1,n-1) = \frac{(-1)^{n-1}}{\prod_{j=0}^{n-1} x_j} g_{|s|}\left(\frac{1}{x_0}, \frac{1}{x_1}, ..., \frac{1}{x_{n-1}}\right) \tag{9}$$

holds. Then from (2), (7) and (9) it follows that

$$f(s-1,n) = \frac{1}{x_n}[f(s,n) - f(s-1,n-1)] =$$

$$= \frac{1}{x_n}\left[\frac{(-1)^n}{\prod_{j=0}^{n} x_j} g_{|s|-1}\left(\frac{1}{x_0}, \frac{1}{x_1}, ..., \frac{1}{x_n}\right) - \frac{(-1)^{n-1}}{\prod_{j=0}^{n-1} x_j} g_{|s|}\left(\frac{1}{x_0}, \frac{1}{x_1}, ..., \frac{1}{x_{n-1}}\right)\right] =$$

$$= \frac{(-1)^n}{x_n \prod_{j=0}^{n-1} x_j}\left[\frac{1}{x_n} g_{|s|-1}\left(\frac{1}{x_0}, \frac{1}{x_1}, ..., \frac{1}{x_n}\right) + g_{|s|}\left(\frac{1}{x_0}, \frac{1}{x_1}, ..., \frac{1}{x_{n-1}}\right)\right] =$$

$$= \frac{(-1)^n}{\prod_{j=0}^{n} x_j} g_{|s|}\left(\frac{1}{x_0}, \frac{1}{x_1}, ..., \frac{1}{x_n}\right). \tag{10}$$

Thus, using mathematical induction from condition (8), assumption (9) and equality (10) proved on the basis of (9), we obtain the identity

$$f(s-1,n) = \frac{(-1)^n}{\prod_{j=0}^{n} x_j} g_{|s|}\left(\frac{1}{x_0}, \frac{1}{x_1}, ..., \frac{1}{x_n}\right), \ n \in \mathbb{N}. \tag{11}$$

Also, using mathematical induction from condition (6), assumption (7) and equality (11) proved on the basis of (7), we obtain the identity

$$f(s,n) = \frac{(-1)^n}{\prod_{j=0}^{n} x_j} g_{|s|-1}\left(\frac{1}{x_0}, \frac{1}{x_1}, ..., \frac{1}{x_n}\right), \ s \leq -1, \ n \in \mathbb{N}.$$

In the case $s \geq 0$, using mathematical induction, using equality (2) and the condition $f(0,n) = 0$, $n \in \mathbb{N}$, (see [11], Appendix E), one can easily prove the identity

$$f(s,n) = 0, \ s \geq 0, \ n \geq s+1. \tag{12}$$

It can also be easily shown that

$$f(s,1) = g_{s-1}(x_0, x_1), \ s \geq 1. \tag{13}$$

From (2), (12) and (13) with $s = n$ we obtain the equality

$$f(n,n) = f(n-1,n-1) + f(n-1,n)x_n = f(n-1,n-1) = ... = f(1,1) =$$

$$= g_0(x_0, x_1) = 1 = g_0(x_0, x_1, ..., x_n). \tag{14}$$

Let us assume that for $s = n+k-1$, where $k \in \mathbb{N}$ is some fixed number, the identity

$$f(n+k-1,n) = g_{k-1}(x_0, x_1, ..., x_n), \ n \in \mathbb{N}, \tag{15}$$

holds. Note that from (13) follows the equality
$$f(1+k,1) = g_k(x_0, x_1). \tag{16}$$
Now let us assume that for some $n \geq 2$ and $s = n-1+k$ the equality
$$f(n-1+k, n-1) = g_k(x_0, x_1, ..., x_{n-1}) \tag{17}$$
holds. Then from (2), (15) and (17) we obtain
$$f(n+k,n) = f(n+k-1, n-1) + f(n+k-1, n)x_n =$$
$$= g_k(x_0, x_1, ..., x_{n-1}) + g_{k-1}(x_0, x_1, ..., x_n)x_n = g_k(x_0, x_1, ..., x_n). \tag{18}$$
Thus, using mathematical induction from condition (16), assumption (17) and equality (18) proved on the basis of (17), we obtain the identity
$$f(n+k,n) = g_k(x_0, x_1, ..., x_n), \ n \in \mathbb{N}. \tag{19}$$
Also, using mathematical induction from condition (14), assumption (15) and equality (19) proved on the basis of (15), we obtain the identity
$$f(n+k,n) = g_k(x_0, x_1, ..., x_n), \ k \in \mathbb{N} \cup \{0\}, \ n \in \mathbb{N}.$$
The lemma is proved by the first method.

The second method. Let $0 \leq s \leq n+1$. The function
$$h_1(x) \stackrel{df}{=} \frac{x^s}{(x-x_0)(x-x_1)...(x-x_n)} - \sum_{j=0}^{n} \frac{x_j^s}{(x-x_j)\prod_{\substack{i=0 \\ i \neq j}}^{n}(x_j - x_i)}, \ x \in \mathbb{C},$$
is analytic and bounded in the entire complex plane, so by Liouville's theorem for functions of a complex variable, it is constant. As $x \to \infty$, we find
$$h_1(x) = \begin{cases} 0, & 0 \leq s \leq n, \\ 1, & s = n+1. \end{cases}$$
From this follows the equality
$$h_1(x_{n+1}) = \sum_{j=0}^{n+1} \frac{x_j^s}{\prod_{\substack{i=0 \\ i \neq j}}^{n+1}(x_j - x_i)} = \begin{cases} 0, & 0 \leq s \leq n, \\ 1, & s = n+1. \end{cases}$$

For $s > n+1$ we get
$$\frac{x^s}{(x-x_0)(x-x_1)...(x-x_n)} = \frac{x^{s-(n+1)}}{\left(1-\frac{x_0}{x}\right)\left(1-\frac{x_1}{x}\right)...\left(1-\frac{x_n}{x}\right)} =$$
$$= x^{s-(n+1)} \left(1 + \frac{x_0}{x} + \frac{x_0^2}{x^2} + \frac{x_0^3}{x^3} + ...\right) \cdot ... \cdot \left(1 + \frac{x_n}{x} + \frac{x_n^2}{x^2} + \frac{x_n^3}{x^3} + ...\right) =$$
$$= x^{s-(n+1)} \left(1 + \frac{g_1(x_0, ..., x_n)}{x} + \frac{g_2(x_0, ..., x_n)}{x^2} + \frac{g_3(x_0, ..., x_n)}{x^3} + ...\right) =$$
$$= x^{s-(n+1)} + x^{s-(n+1)-1} g_1(x_0, ..., x_n) + x^{s-(n+1)-2} g_2(x_0, ..., x_n) + ... +$$
$$+ g_{s-(n+1)}(x_0, ..., x_n) + O\left(\frac{1}{x}\right), \ x \to \infty.$$
Therefore, the function
$$h_2(x) \stackrel{df}{=} \frac{x^s}{(x-x_0)(x-x_1)...(x-x_n)} - \sum_{j=0}^{n} \frac{x_j^s}{(x-x_j)\prod_{\substack{i=0 \\ i \neq j}}^{n}(x_j - x_i)} -$$
$$- \Big[ x^{s-(n+1)} + x^{s-(n+1)-1} g_1(x_0, ..., x_n) + x^{s-(n+1)-2} g_2(x_0, ..., x_n) + ... +$$

$$+ g_{s-(n+1)}(x_0,...,x_n)\Big]$$

is analytic in the entire complex plane and tends to zero as $x \to \infty$. Consequently, by Liouville's theorem for functions of a complex variable, the identity $h_2(x) \equiv 0$ holds, in particular, $h_2(x_{n+1}) = 0$. From the last equality, we obtain

$$\sum_{j=0}^{n+1} \frac{x_j^s}{\prod_{\substack{i=0 \\ i \neq j}}^{n+1}(x_j - x_i)} = x_{n+1}^{s-(n+1)} + x_{n+1}^{s-(n+1)-1} g_1(x_0,...,x_n) + x_{n+1}^{s-(n+1)-2} g_2(x_0,...,x_n) + ... +$$

$$+ g_{s-(n+1)}(x_0,...,x_n) = g_{s-(n+1)}(x_0,...,x_{n+1}).$$

Let $s < 0$, then

$$\frac{x^s}{(x-x_0)(x-x_1)...(x-x_n)} = \frac{1}{x^{|s|}(-1)^{n+1}\left(\prod_{k=0}^{n}x_k\right)\left(1-\frac{x}{x_0}\right)\left(1-\frac{x}{x_1}\right)...\left(1-\frac{x}{x_n}\right)} =$$

$$= \frac{(-1)^{n+1}}{\left(\prod_{k=0}^{n}x_k\right)x^{|s|}}\left(1+\frac{x}{x_0}+\frac{x^2}{x_0^2}+\frac{x^3}{x_0^3}+...\right)\cdot...\cdot\left(1+\frac{x}{x_n}+\frac{x^2}{x_n^2}+\frac{x^3}{x_n^3}+...\right) =$$

$$= \frac{(-1)^{n+1}}{\left(\prod_{k=0}^{n}x_k\right)x^{|s|}}\left(1+xg_1\left(\frac{1}{x_0},...,\frac{1}{x_n}\right)+x^2 g_2\left(\frac{1}{x_0},...,\frac{1}{x_n}\right)+x^3 g_3\left(\frac{1}{x_0},...,\frac{1}{x_n}\right)+...\right) =$$

$$= \frac{(-1)^{n+1}}{\left(\prod_{k=0}^{n}x_k\right)}\Bigg[\frac{1}{x^{|s|}} + \frac{1}{x^{|s|-1}}g_1\left(\frac{1}{x_0},...,\frac{1}{x_n}\right) + \frac{1}{x^{|s|-2}}g_2\left(\frac{1}{x_0},...,\frac{1}{x_n}\right)+...+$$

$$+\frac{1}{x}g_{|s|-1}\left(\frac{1}{x_0},...,\frac{1}{x_n}\right) + O(1)\Bigg], \ x \to 0.$$

Therefore, the function

$$h_3(x) \stackrel{df}{=} \frac{x^s}{(x-x_0)(x-x_1)...(x-x_n)} - \sum_{j=0}^{n}\frac{x_j^s}{(x-x_j)\prod_{\substack{i=0 \\ i \neq j}}^{n}(x_j - x_i)} -$$

$$-\frac{(-1)^{n+1}}{\left(\prod_{k=0}^{n}x_k\right)x}\Bigg[\frac{1}{x^{|s|-1}} + \frac{1}{x^{|s|-2}}g_1\left(\frac{1}{x_0},...,\frac{1}{x_n}\right) + \frac{1}{x^{|s|-3}}g_2\left(\frac{1}{x_0},...,\frac{1}{x_n}\right)+...+$$

$$+g_{|s|-1}\left(\frac{1}{x_0},...,\frac{1}{x_n}\right)\Bigg]$$

is analytic in the entire complex plane and tends to zero as $x \to \infty$. Then, by Liouville's theorem for functions of a complex variable, the identity $h_3(x) \equiv 0$ holds, in particular, $h_3(x_{n+1}) = 0$. From the last equality, we obtain

$$\sum_{j=0}^{n+1} \frac{x_j^s}{\prod_{\substack{i=0 \\ i \neq j}}^{n+1}(x_j - x_i)} = \frac{(-1)^{n+1}}{\prod_{k=0}^{n+1}x_k}\Bigg[\frac{1}{x_{n+1}^{|s|-1}} + \frac{1}{x_{n+1}^{|s|-2}}g_1\left(\frac{1}{x_0},...,\frac{1}{x_n}\right)+$$

$$+\frac{1}{x_{n+1}^{|s|-3}} g_2\left(\frac{1}{x_0},...,\frac{1}{x_n}\right)+...+ g_{|s|-1}\left(\frac{1}{x_0},...,\frac{1}{x_n}\right)\right] = \frac{(-1)^{n+1}}{\prod_{k=0}^{n+1} x_k} g_{|s|-1}\left(\frac{1}{x_0},...,\frac{1}{x_{n+1}}\right).$$

In the previous discussion, we considered variables $x_0$, $x_1$, ..., $x_{n+1}$ for $n \in \mathbb{N}$. The case of two variables $x_0$, $x_1$ is elementary.

The lemma is proved using the second method. The lemma is proven.

**3. Volterra's periodic chain.** Assume that $u_n > 0$, $n \in \mathbb{Z}$, and make the change of variables [4.1] from [6] in system (1)

$$a_n = \frac{1}{2}\sqrt{u_n},$$

then

$$a'_n = 2a_n\left(a_{n+1}^2 - a_{n-1}^2\right), \; n \in \mathbb{Z}. \tag{20}$$

In what follows, the theory of the discrete Hill equation [9-11] will be used substantially in the notation adopted in [9]. A detailed exposition of this extensive theory can be found in the cited works, and therefore, for brevity, it is not presented here. We also note that this note can be considered as an adaptation to the Volterra chain of the approach proposed in [9] for studying the Toda chain.

In the case of the Volterra chain in the discrete Hill equation (4) from [9], the coefficients $a_n$ are solutions of system (20) that satisfy the periodicity condition $a_{n+N} = a_n$, where $N \in \mathbb{N}$ is the period of system (1), $n \in \mathbb{Z}$; the coefficients $b_n = 0$, $n \in \mathbb{Z}$. For the spectral parameters $\mu_j$, $j = \overline{1, N-1}$, of the Hill equation (4) from [9], the differential identity (24) in [9] is satisfied; we write it out and transform it taking into account the above for the coefficients $a_n$, $b_n$:

$$\mu'_j(t) = \sum_{n=1}^{N}\left(2a'_n(t)y_n^j y_{n+1}^j + b'_n(t)\left(y_n^j\right)^2\right) = \sum_{n=1}^{N} 2a'_n(t)y_n^j y_{n+1}^j =$$

$$= \sum_{n=1}^{N} 4a_n\left(a_{n+1}^2 - a_{n-1}^2\right)y_n^j y_{n+1}^j = \sum_{n=1}^{N} 4a_n a_{n+1}^2 y_n^j y_{n+1}^j - \sum_{n=1}^{N} 4a_n a_{n-1}^2 y_n^j y_{n+1}^j =$$

$$= \sum_{n=1}^{N} 4a_n a_{n+1} y_n^j \left(a_{n+1} y_{n+1}^j\right) - \sum_{n=1}^{N} 4a_n a_{n-1}\left(a_{n-1} y_n^j\right) y_{n+1}^j =$$

$$= \sum_{n=1}^{N} 4a_n a_{n+1} y_n^j \left(\mu_j(t) y_{n+2}^j - a_{n+2} y_{n+3}^j\right) - \sum_{n=1}^{N} 4a_n a_{n-1}\left(\mu_j(t) y_{n-1}^j - a_{n-2} y_{n-2}^j\right) y_{n+1}^j =$$

$$= \sum_{n=1}^{N} 4\mu_j(t) a_n a_{n+1} y_n^j y_{n+2}^j - \sum_{n=1}^{N} 4a_n a_{n+1} a_{n+2} y_n^j y_{n+3}^j -$$

$$-\sum_{n=1}^{N} 4\mu_j(t) a_{n-1} a_n y_{n-1}^j y_{n+1}^j + \sum_{n=1}^{N} 4a_{n-2} a_{n-1} a_n y_{n-2}^j y_{n+1}^j =$$

$$= \sum_{n=1}^{N} 4\mu_j(t) a_n a_{n+1} y_n^j y_{n+2}^j - \sum_{n=1}^{N} 4a_n a_{n+1} a_{n+2} y_n^j y_{n+3}^j -$$

$$-\sum_{n=0}^{N-1} 4\mu_j(t) a_n a_{n+1} y_n^j y_{n+2}^j + \sum_{n=-1}^{N-2} 4a_n a_{n+1} a_{n+2} y_n^j y_{n+3}^j =$$

$$= 4\mu_j(t) a_N a_{N+1} y_N^j y_{N+2}^j - \sum_{n=N-1}^{N} 4a_n a_{n+1} a_{n+2} y_n^j y_{n+3}^j -$$

$$-4\mu_j(t) a_0 a_1 y_0^j y_2^j + \sum_{n=-1}^{0} 4a_n a_{n+1} a_{n+2} y_n^j y_{n+3}^j. \tag{21}$$

From the Hill equation (4) in [9] and the condition $y_1^j = y_{N+1}^j = 0$ in identity (24) from [9] we obtain the equalities $a_{-1} y_{-1}^j = \mu_j y_0^j$ and $a_{N-1} y_{N-1}^j = \mu_j y_N^j$. Therefore,

$$4a_na_{n+1}a_{n+2}y_n^j y_{n+3}^j \big|_{n=-1} = 4a_{-1}a_0 a_1 y_{-1}^j y_2^j = 4\mu_j(t)a_0 a_1 y_0^j y_2^j,$$

$$4a_na_{n+1}a_{n+2}y_n^j y_{n+3}^j \big|_{n=N-1} = 4a_{N-1}a_N a_{N+1} y_{N-1}^j y_{N+2}^j = 4\mu_j(t)a_N a_{N+1} y_N^j y_{N+2}^j.$$

Again from Hill's equation (4) in [9] and the condition $y_1^j = y_{N+1}^j = 0$ in identity (24) from [9] we obtain the equalities

$$a_1 y_1^j + 0 \cdot y_2^j + a_2 y_3^j = \mu_j y_2^j, \ a_2 y_3^j = \mu_j y_2^j;$$
$$a_0 y_0^j + 0 \cdot y_1^j + a_1 y_2^j = \mu_j y_1^j, \ a_1 y_2^j = -a_0 y_0^j.$$

Then

$$4a_n a_{n+1} a_{n+2} y_n^j y_{n+3}^j \big|_{n=0} = 4a_0 a_1 a_2 y_0^j y_3^j = 4\mu_j(t) a_0 a_1 y_0^j y_2^j = -4\mu_j(t) a_0^2 \left(y_0^j\right)^2.$$

Similarly, from the Hill equation (4) in [9] and the condition $y_1^j = y_{N+1}^j = 0$ in identity (24) from [9] we obtain the equalities

$$a_{N+1} y_{N+1}^j + 0 \cdot y_{N+2}^j + a_{N+2} y_{N+3}^j = \mu_j y_{N+2}^j, \ a_{N+2} y_{N+3}^j = \mu_j y_{N+2}^j;$$
$$a_N y_N^j + 0 \cdot y_{N+1}^j + a_{N+1} y_{N+2}^j = \mu_j y_{N+1}^j, \ a_{N+1} y_{N+2}^j = -a_N y_N^j.$$

Therefore

$$4a_n a_{n+1} a_{n+2} y_n^j y_{n+3}^j \big|_{n=N} = 4a_N a_{N+1} a_{N+2} y_N^j y_{N+3}^j = 4\mu_j(t) a_N a_{N+1} y_N^j y_{N+2}^j =$$
$$= -4\mu_j(t) a_N^2 \left(y_N^j\right)^2.$$

Taking into account the equalities obtained in the previous paragraph, we continue calculating the derivative $\mu_j'(t)$ in (21):

$$\mu_j'(t) = 4\mu_j(t) a_0^2 \left[\left(y_N^j\right)^2 - \left(y_0^j\right)^2\right]. \tag{22}$$

Let us write down the equalities on page 22 in [9], which we will need later:

$$\|\theta^j\|^2 = \sum_{n=1}^N \left(\theta_n^j\right)^2 = a_N \theta_N^j \theta_{N+1}'\big|_{\lambda=\mu_j}, \ \theta' = \frac{d\theta}{d\lambda},$$

$$\left(y_0^j\right)^2 = \frac{\left(\theta_0^j\right)^2}{\|\theta^j\|^2}, \ \left(y_N^j\right)^2 = \frac{\left(\theta_N^j\right)^2}{\|\theta^j\|^2};$$

note that the first of these equalities is obtained, in particular, from equalities [3.2] and [3.3] in [6], in this case, in equality [3.2] from [6] it is necessary to put the fraction $\dfrac{a_1}{a_0}$ instead of the fraction $\dfrac{a_0}{a_1}$ (see also (4.1.38) in [11]).

Taking into account the last four equalities, the initial condition $\theta_0^j = 1$, and the periodicity of the coefficients $a_n$, from (22) we obtain

$$\mu_j'(t) = 4\mu_j(t) a_0^2 \left[\frac{\left(\theta_N^j\right)^2}{\|\theta^j\|^2} - \frac{\left(\theta_0^j\right)^2}{\|\theta^j\|^2}\right] = 4\mu_j(t) a_0 \left[\theta_N^j - \frac{1}{\theta_N^j}\right] \frac{1}{\theta_{N+1}'\big|_{\lambda=\mu_j}}. \tag{23}$$

Using the equality

$$\theta_N(\lambda,t)\varphi_{N+1}(\lambda,t) - \theta_{N+1}(\lambda,t)\varphi_N(\lambda,t) = 1$$

on page 22 in [9] (see also [2.10] in [6] or (4.1.15) in [11]), we obtain the following identities for the function $\Delta$ from [9] (see also [2.8], [2.9] in [6] or (4.1.19) in [11]):

$$\Delta^2(\mu_j(t)) - 4 = \left[\theta_N(\mu_j(t),t) - \varphi_{N+1}(\mu_j(t),t)\right]^2 +$$
$$+ 4\theta_N(\mu_j(t),t)\varphi_{N+1}(\mu_j(t),t) - 4 =$$
$$= \left[\theta_N(\mu_j(t),t) - \varphi_{N+1}(\mu_j(t),t)\right]^2 = \left[\theta_N(\mu_j(t),t) - \frac{1}{\theta_N(\mu_j(t),t)}\right]^2,$$

$$\theta_N(\mu_j(t),t) - \frac{1}{\theta_N(\mu_j(t),t)} = \sigma_j(t)\sqrt{\Delta^2(\mu_j(t))-4},$$

where $\sigma_j(t) = \text{sign}\left(\theta_N(\mu_j(t),t) - \frac{1}{\theta_N(\mu_j(t),t)}\right)$, $j = \overline{1, N-1}$. Using formula (7) in [9], from the last equality we obtain

$$\theta_N^j(\mu_j(t),t) - \frac{1}{\theta_N^j(\mu_j(t),t)} = \sigma_j(t)\left(\prod_{k=1}^N a_k\right)^{-1}\sqrt{\prod_{i=1}^{2N}(\mu_j - \lambda_i)}.$$

From formula (8) in [9] follows the identity

$$\theta'_{N+1}(\lambda)\big|_{\lambda=\mu_j} = -a_0\left(\prod_{k=1}^N a_k\right)^{-1}\prod_{\substack{i=1\\i\neq j}}^{N-1}(\mu_j - \mu_i).$$

Substituting the last two relations into (23), we obtain

$$\mu'_j(t) = -4\mu_j(t)\sigma_j(t)\frac{\sqrt{\prod_{i=1}^{2N}(\mu_j - \lambda_i)}}{\prod_{\substack{i=1\\i\neq j}}^{N-1}(\mu_j - \mu_i)}, \quad j = \overline{1, N-1}. \tag{24}$$

Let us show that the numbers $\lambda_j(t)$, $j = \overline{1, 2N}$, for the discrete Hill equation (4) from [9] do not depend on $t$. From equality (33) in [9] it follows that for these numbers the differential identity

$$\lambda'_j(t) = \sum_{n=1}^N \left(2a'_n(t)y_n^j y_{n+1}^j + b'_n(t)\left(y_n^j\right)^2\right)$$

holds, in which some non-zero sequence $y_n^j \in \mathbb{C}$, $n \in \mathbb{Z}$, is periodic $y_n^j = y_{n+N}^j$ or antiperiodic $y_n^j = -y_{n+N}^j$; we transform it taking into account the above for the quantities $y_n^j$ and the coefficients $a_n$, $b_n$:

$$\lambda'_j(t) = \sum_{n=1}^N \left(2a'_n(t)y_n^j y_{n+1}^j + b'_n(t)\left(y_n^j\right)^2\right) = \sum_{n=1}^N 2a'_n(t)y_n^j y_{n+1}^j =$$

$$= \sum_{n=1}^N 4a_n\left(a_{n+1}^2 - a_{n-1}^2\right)y_n^j y_{n+1}^j = \sum_{n=1}^N 4a_n a_{n+1}^2 y_n^j y_{n+1}^j - \sum_{n=1}^N 4a_n a_{n-1}^2 y_n^j y_{n+1}^j =$$

$$= \sum_{n=1}^N 4a_n a_{n+1} y_n^j \left(a_{n+1} y_{n+1}^j\right) - \sum_{n=1}^N 4a_n a_{n-1}\left(a_{n-1} y_n^j\right) y_{n+1}^j =$$

$$= \sum_{n=1}^N 4a_n a_{n+1} y_n^j \left(\lambda_j(t) y_{n+2}^j - a_{n+2} y_{n+3}^j\right) - \sum_{n=1}^N 4a_n a_{n-1}\left(\lambda_j(t) y_{n-1}^j - a_{n-2} y_{n-2}^j\right) y_{n+1}^j =$$

$$= \sum_{n=1}^N 4\lambda_j(t) a_n a_{n+1} y_n^j y_{n+2}^j - \sum_{n=1}^N 4a_n a_{n+1} a_{n+2} y_n^j y_{n+3}^j -$$

$$- \sum_{n=1}^N 4\lambda_j(t) a_{n-1} a_n y_{n-1}^j y_{n+1}^j + \sum_{n=1}^N 4a_{n-2} a_{n-1} a_n y_{n-2}^j y_{n+1}^j =$$

$$= \sum_{n=1}^N 4\lambda_j(t) a_n a_{n+1} y_n^j y_{n+2}^j - \sum_{n=1}^N 4a_n a_{n+1} a_{n+2} y_n^j y_{n+3}^j -$$

$$- \sum_{n=0}^{N-1} 4\lambda_j(t) a_n a_{n+1} y_n^j y_{n+2}^j + \sum_{n=-1}^{N-2} 4a_n a_{n+1} a_{n+2} y_n^j y_{n+3}^j =$$

$$= 4\lambda_j(t) a_N a_{N+1} y_N^j y_{N+2}^j - \sum_{n=N-1}^N 4a_n a_{n+1} a_{n+2} y_n^j y_{n+3}^j -$$

$$-4\lambda_j(t)a_0 a_1 y_0^j y_2^j + \sum_{n=-1}^{0} 4a_n a_{n+1} a_{n+2} y_n^j y_{n+3}^j = 0.$$

From equality (12) in [9] and the condition $b_n = 0$, $n \in \mathbb{Z}$, for the Volterra chain we obtain the identities

$$\sum_{j=1}^{N-1} \mu_j^{(n)} = 0, \quad n \geq 0. \tag{25}$$

In particular, from the last equality for $n = 1$ and equations (24) follows the identity

$$\sum_{j=1}^{N-1} \mu_j(t)\sigma_j(t) \frac{\sqrt{\prod_{i=1}^{2N}(\mu_j - \lambda_i)}}{\prod_{\substack{i=1\\i \neq j}}^{N-1}(\mu_j - \mu_i)} = 0. \tag{26}$$

Note that the sum $\sum_{j=1}^{N-1} \mu_j$ on the left-hand side of (25) for $n = 0$ and the left-hand side of the last equality are not first integrals for system (24). In the case of $N = 3$, equality (26) can be easily verified if we take into account that $\mu_1 + \mu_2 = 0$ (see (25) for $n = 0$), $\Delta(\lambda) = \frac{1}{a_1 a_2 a_3}\left[\lambda^3 - \left(\sum_{j=1}^{3} a_j^2\right)\lambda\right]$ and therefore

$\{\lambda_1, ..., \lambda_6\} = \{-\lambda_1, ..., -\lambda_6\}$, $\theta_3(\lambda) = -\frac{a_0}{a_1 a_2}\lambda$ and, consequently, $\sigma_1 = -\sigma_2$.

Let $N \geq 4$. Then, using the lemma and condition (25) for $n = 0$, we write the following equalities

$$\sum_{j=1}^{N-1} \frac{\mu_j^s}{\prod_{\substack{i=1\\i \neq j}}^{N-1}(\mu_j - \mu_i)} = 0, \quad s = \overline{1, N-3}; \quad \sum_{j=1}^{N-1} \frac{\mu_j^{N-2}}{\prod_{\substack{i=1\\i \neq j}}^{N-1}(\mu_j - \mu_i)} = 1,$$

$$\sum_{j=1}^{N-1} \frac{\mu_j^{N-1}}{\prod_{\substack{i=1\\i \neq j}}^{N-1}(\mu_j - \mu_i)} = \sum_{j=1}^{N-1} \mu_j = 0.$$

Using equations (24), the last identities can be rewritten as

$$\sum_{j=1}^{N-1} \frac{\mu_j^{s-1}\mu_j'(t)}{\sigma_j(t)\sqrt{\prod_{i=1}^{2N}(\mu_j - \lambda_i)}} = 0, \quad s = \overline{1, N-3}; \quad \sum_{j=1}^{N-1} \frac{\mu_j^{N-3}\mu_j'(t)}{\sigma_j(t)\sqrt{\prod_{i=1}^{2N}(\mu_j - \lambda_i)}} = -4,$$

$$\sum_{j=1}^{N-1} \frac{\mu_j^{N-2}\mu_j'(t)}{\sigma_j(t)\sqrt{\prod_{i=1}^{2N}(\mu_j - \lambda_i)}} = 0. \tag{27}$$

System (27) is close in form to system [4.1a], [4.1b] in [12], which can be solved using the classical theory of hyperelliptic (Abelian) integrals.

Let the period in system (1) be odd, $N = 2l + 1$, $l \geq 0$. Then, following the approach proposed in [10], from equality (15) in [9] we obtain the identities

$$a_k^2 + a_{k+1}^2 = \frac{1}{4}\sum_{j=1}^{2N}\lambda_j^2 - \frac{1}{2}\sum_{j=1}^{N-1}\mu_{j,k}^2, \quad k \in \mathbb{Z};$$

$$a_m^2 = \frac{1}{2}\sum_{k=0}^{N-1}\left[(-1)^k\left(\frac{1}{4}\sum_{j=1}^{2N}\lambda_j^2 - \frac{1}{2}\sum_{j=1}^{N-1}\mu_{j,m+k}^2\right)\right] =$$

$$= \frac{1}{8}\left(\sum_{j=1}^{2N}\lambda_j^2\right) - \frac{1}{4}\sum_{k=0}^{N-1}\left[(-1)^k\sum_{j=1}^{N-1}\mu_{j,m+k}^2\right], \quad m \in \mathbb{Z}, \tag{28}$$

where the quantities $\mu_{j,k}$, $j = \overline{1, N-1}$, are determined by equation (19) in [9].

In the case of an arbitrary period $N \in \mathbb{N}$ in system (1), from equalities (20), (21) in [9] and the condition $b_n = 0$, $n \in \mathbb{Z}$, for the Volterra chain, we obtain the identities

$$a_k^2 = \frac{1}{8}\sum_{j=1}^{2N} \lambda_j^2 - \frac{1}{4}\sum_{j=1}^{N-1}\mu_{j,k}^2 - \frac{1}{2}\sum_{j=1}^{N-1}\frac{\sigma_{j,k}\sqrt{\prod_{i=1}^{2N}(\mu_{j,k} - \lambda_i)}}{\prod_{\substack{i=1\\i\neq j}}^{N-1}(\mu_{j,k} - \mu_{i,k})}, \quad k \in \mathbb{Z}. \tag{29}$$

As for the spectral parameters $\mu_j$, we similarly calculate the differential equations for the quantities $\mu_{j,k}$:

$$\mu'_{j,k}(t) = -4\mu_{j,k}(t)\sigma_{j,k}(t)\frac{\sqrt{\prod_{i=1}^{2N}(\mu_{j,k} - \lambda_i)}}{\prod_{\substack{i=1\\i\neq j}}^{N-1}(\mu_{j,k} - \mu_{i,k})}, \quad j = \overline{1, N-1},\ k \in \mathbb{Z}, \tag{30}$$

where $\sigma_{j,k}(t) = sign\left(\theta_{N,k}(\mu_{j,k}(t), t) - \frac{1}{\theta_{N,k}(\mu_{j,k}(t), t)}\right)$, the functions $\theta_{N,k}$ are determined by equation (19) in [9]. Using equalities (30), we write identities (29) in the following form:

$$a_k^2 = \frac{1}{8}\sum_{j=1}^{2N}\lambda_j^2 - \frac{1}{4}\sum_{j=1}^{N-1}\mu_{j,k}^2 + \frac{1}{8}\frac{d}{dt}\ln\left(\prod_{j=1}^{N-1}\mu_{j,k}(t)\right). \tag{31}$$

The method for calculating symmetric polynomials of variables $\mu_{j,k}$, $j = \overline{1, N-1}$, in formulas (28) and (31) is proposed in [10, 11] (see also [8]).

This work was supported by a grant from the Simons Foundation (SFI-PD-Ukraine-00014586, B.D.V.).